\documentclass[12pt]{article}
\usepackage{a4}
\usepackage{amsfonts}
\usepackage{amssymb}
\usepackage{color}

\begin{document}
\title{\textbf{On the structure of the fibers of truncation morphisms}}
\author{Helena Cobo\thanks{Supported by FWO-Flanders project GO31806N and
by MCI grant MTM2010-21740\mbox{-}C02\mbox{-}01} \and Dirk
Segers\thanks{Postdoctoral Fellow of the Fund for Scientific
Research - Flanders (Belgium).
}}

\date{August 11, 2012}

\maketitle

\begin{abstract}
Let $k$ be an algebraically closed field and let $X$ be a
separated scheme of finite type over $k$ of pure dimension $d$. We
study the structure of the fibres of the truncation morphisms from
the arc space of $X$ to jet spaces of $X$ and also between jet
spaces. Our results are generalizations of results of Denef,
Loeser, Ein and Musta\c t\v a. We will use them to find the
optimal lower bound for the poles of the motivic zeta function
associated to an arbitrary ideal.
\end{abstract}

\section{Introduction}

Let $k$ be an algebraically closed field of arbitrary
characteristic. The functor $\cdot \times_{\mbox{\scriptsize{Spec
}} k} \mbox{Spec } k[t]/(t^{n+1})$ on the category of schemes over
$k$  has a right adjoint, which we denote by $\mathcal{L}_n$. Let
$X$ be a separated scheme of finite type over $k$. We call
$\mathcal{L}_n(X)$ the scheme of $n$-jets of $X$ and we define an
$n$-jet on $X$ as a closed point on $\mathcal{L}_n(X)$. An $n$-jet
on $X$ is thus a morphism $\mbox{Spec}(k[t]/(t^{n+1})) \rightarrow
X$. For $m,n \in \mathbb{Z}_{\geq 0}$ satisfying $m \geq n$, the
canonical embeddings $X \times_{\mbox{\scriptsize{Spec }} k}
\mbox{Spec } k[t]/(t^{n+1}) \hookrightarrow X
\times_{\mbox{\scriptsize{Spec }} k} \mbox{Spec } k[t]/(t^{m+1})$
induce canonical truncation morphisms $\pi_n^m : \mathcal{L}_m(X)
\rightarrow \mathcal{L}_n(X)$. The projective limit of the
projective system
\[
\cdots \rightarrow \mathcal{L}_{n+1}(X) \rightarrow \mathcal{L}_n(X)
\rightarrow \cdots \rightarrow \mathcal{L}_0(X)=X
\]
exists in the category of schemes over $k$. It is denoted by
$\mathcal{L}(X)$ and it is called the arc space of $X$. An arc on
$X$ is a closed point on $\mathcal{L}(X)$, so it is a morphism
$\mbox{Spec}(k[[t]]) \rightarrow X$. In general, $\mathcal{L}(X)$
is not of finite type over $k$. There are natural truncation
morphisms $\pi_n : \mathcal{L}(X) \rightarrow \mathcal{L}_n(X)$.
For more information about these constructions, see
\cite{DenefLoeserGermsofarcs} or \cite{Mustata}.

In this paper we study some properties of the truncation
morphisms. If $X$ is a nonsingular variety of dimension $d$, then
$\pi_n(\mathcal{L}(X)) = \mathcal{L}_n(X)$ and the truncation
morphisms $\pi^{n+1}_n : \mathcal{L}_{n+1}(X) \rightarrow
\mathcal{L}_n(X)$ are locally trivial fibrations with fibre
$\mathbb{A}^d$. The general case is more complicated, and partial
results generalizing the nonsingular case are already in
\cite{DenefLoeserGermsofarcs}, \cite{Looijenga} and
\cite{EinMustata}. We first describe these results and thereafter we formulate our generalizations.

The first result concerns the truncation map $\pi_n (\mathcal{L}(X)) \rightarrow \pi_l(\mathcal{L}(X))$ for $n>l$ and is due to Denef and Loeser. Let $X$ be a separated scheme of finite type over $k$ of pure dimension $d$. Let $\mbox{Jac}_X$ be the Jacobian subscheme of $X$ and for a positive integer $e$, we denote by $\mbox{Cont}^e(\mbox{Jac}_X)$ the set of all $\gamma \in \mathcal{L}(X)$ satisfying $\mbox{ord}_{\gamma}(\mbox{Jac}_X)=e$. It is proved in \cite{DenefLoeserGermsofarcs} that for $n > l \geq e$, the truncation map
\[
\pi_n(\mbox{Cont}^e(\mbox{Jac}_X)) \rightarrow
\pi_l(\mbox{Cont}^e(\mbox{Jac}_X))
\]
is a piecewise trivial fibration with fibre $\mathbb{A}^{(n-l)d}$.

The second result concerns the truncation map $\mathcal{L}_n(X) \rightarrow \pi^n_l(\mathcal{L}_n(X))$ for $n>l$ and is due to Ein and Musta\c t\v a. Let $X$ be a separated scheme of finite type over $k$ of pure dimension $d$ which is either reduced or a locally complete intersection. For positive integers $n,e$ with $n \geq e$, we denote by $\mbox{Cont}^e(\mbox{Jac}_X)_n$ the set of all $\theta \in \mathcal{L}_n(X)$ satisfying $\mbox{ord}_{\theta}(\mbox{Jac}_X)=e$. Note that this is well defined because $n \geq e$. It is proved in \cite{EinMustata} that for positive integers $n$, $l$ and $e$ satisfying $l+e \leq n \leq 2l$, the truncation map
\[
\mbox{Cont}^e(\mbox{Jac}_X)_n \rightarrow
\pi^n_l(\mbox{Cont}^e(\mbox{Jac}_X)_n)
\]
is a piecewise trivial fibration with fibre
$\mathbb{A}^{(n-l)d+e}$. 

The above results involve the invariant $e$ of an arc or jet with generic point on the smooth part of $X$. To describe our results, we need to refine this invariant. Let $X$ be a closed subscheme of $\mathbb{A}^N$ with ideal generated by $m$ polynomials, say $I_X = (f_1,\ldots,f_m)$. We will suppose in this introduction that $X$ is of pure dimension $d$ and we put $r=N-d$. For an arc $\gamma \in \mathcal{L}(X)$, there exists a unique sequence $0 \leq e_1 \leq e_2 \leq \cdots \leq e_p$ in $\mathbb{N}$, with $0\leq p \leq \min\{m,N\}$, such that $t^{e_1},t^{e_2},\ldots,t^{e_p}$ are the invariant factors of the matrix $((\partial f_i/\partial x_j)(\gamma))_{i,j}$. We have now that $p \leq r$ and equality holds if and only if the generic point of $\gamma$ lies on the smooth part of $X$, in which case $e = e_1 + e_2 + \ldots + e_r$. Note that $p$ and the $e_i$ depend on $X$ and $\gamma$, so we will write $p = p(X;\gamma)$ and $e_i = e_i(X;\gamma)$ if necessary. 
For a sequence $0 \leq e_1 \leq e_2 \leq \cdots \leq e_p$ in $\mathbb{N}$, let $\mathcal{L}^{e_1,\ldots,e_p}(X)$ denote the set of all $\gamma \in \mathcal{L}(X)$ for which $p(X;\gamma) = p$ and $e_i(X;\gamma) = e_i$ for $1 \leq i \leq p$. Note that $\mbox{Cont}^e(\mbox{Jac}_X)$ is partitioned into the sets $\mathcal{L}^{e_1,\ldots,e_r}(X)$, with $e_1+e_2+\ldots+e_r=e$, and this partition is finite. The $n$-jets are treated in a similar way. Let $g$ be $n/2$ if $n$ is even and $(n+1)/2$ if $n$ is odd. To $\theta \in \mathcal{L}_n(X)$, we will associate a sequence $0 \leq e_1 \leq e_2 \leq \cdots \leq e_N \leq n+1$ in $\mathbb{N}$ and a natural number $b(X;\theta)$ which is the largest $i \in \{0,\ldots,N\}$ for which $e_i < g$. We denote $e_i = e_i(X;\theta)$ and we will see in Proposition 6(a) that $b(X;\theta) \leq r$. For a sequence $0 \leq e_1 \leq e_2 \leq \cdots \leq e_b < g$ in $\mathbb{N}$ with $b \leq N$, let $\mathcal{L}^{e_1,\ldots,e_b}_n(X)$ denote the set of all $\theta \in \mathcal{L}_n(X)$ for which $b(X;\theta) = b$ and $e_i(X;\theta) = e_i$ for $1 \leq i \leq b$.

We first state our most important results concerning arcs. For $n > l \geq e_r$ we consider the truncation map $\pi_n(\mathcal{L}^{e_1,\ldots,e_r}(X)) \rightarrow \pi_l(\mathcal{L}^{e_1,\ldots,e_r}(X))$. Proposition 5(b) tells us that this map is a piecewise trivial fibration between constructible sets with fibre $\mathbb{A}^{(n-l)d}$. According to Propositon 1(c), this map is even a locally trivial fibration between locally closed sets in the case $m=r$. Note that our conditions on $n$ and $l$ are less restrictive than in the theorem of Denef and Loeser. We also obtained new sufficient conditions in order to be able to lift a jet to an arc on $X$ (Proposition 1(a) and 7(a)(i)). If $\theta \in \mathcal{L}^{e_1,\ldots,e_r}_n(X)$, we can lift $\pi^{n}_{n-e_r}(\theta)$ to an arc on $X$ if $m=r$ and $n \geq 2e_r$ and also if $X$ is reduced and $n \geq \max \{2e_r,e\}$. The case $m=r$ can be considered as the complete intersection case, and is treated in Section 2. In this section, we obtain results by working with the equations. In section 3, we consider also closed subschemes of $\mathbb{A}^N$ which are not a complete intersection. The proofs in section 3 use intensively the results of section 2, which are applied to complete intersections $M$ given by $r$ well-chosen polynomials in the ideal of $X$. The crucial point is that we can prove that certain arcs constructed on $M$ also lie on $X$.

We obtain similar results concerning jets. For $n \geq 2e_r$ and $n > l \geq e_r$ we consider the truncation map $\mathcal{L}^{e_1,\ldots,e_r}_n(X) \rightarrow \pi^n_l(\mathcal{L}^{e_1,\ldots,e_r}_n(X))$. Proposition 2(c) tells us that this map is a locally trivial fibration between locally closed sets if $m=r$. In Proposition 7(b), we obtain that it is a piecewise trivial fibration if $X$ is reduced and $n \geq e$. In both cases the fibre is an affine space, and if $l \leq n - e_r$, the dimension of this fibre is $(n-l)d+e$. Again, our conditions on $n$ and $l$ are much weaker than in the corresponding theorem of Ein and Musta\c t\v a. We also want to mention that we obtained in Proposition 3(c) and 6(b) similar results for the truncation map $\mathcal{L}^{e_1,\ldots,e_b}_n(X) \rightarrow \pi^n_l(\mathcal{L}^{e_1,\ldots,e_b}_n(X))$ for $n > l \geq n-g$ and $b \leq r$. 

The results in the previous paragraph have important consequences for the motivic zeta function. Here, the Grothendieck ring of algebraic varieties comes into play. The Grothendieck ring $K_0(\mbox{Var}_{k})$ is the abelian group generated by the symbols $[V]$, where $V$ is an algebraic variety, subject to the relations $[V]=[V']$, if $V$ is isomorphic to $V'$, and $[V]=[V \setminus W] + [W]$, if $W$ is closed in $V$. The ring structure of $K_0(\mbox{Var}_{k})$ is given by $[V] \cdot [W] := [V \times W]$. We denote by $\mathbb{L}$ the class of the affine line.

When we consider the class of $\mathcal{L}_n(X)$ in the Grothendieck ring, we endow $\mathcal{L}_n(X)$ with its reduced structure so that it becomes a variety. Our results about the fibres of the truncation morphisms will imply that $[\mathcal{L}_n(X)]$ is a multiple of $\mathbb{L}^{\ulcorner (N-m+1)(n/2) \urcorner}$ in $K_0(\mbox{Var}_{k})$ for all $n \in \mathbb{Z}_{\geq 0}$ if $m<N$ and that it is a multiple of $\mathbb{L}^{\ulcorner (d+1)(n/2) \urcorner}$ in $K_0(\mbox{Var}_{k})$ for all $n \in \mathbb{Z}_{\geq 0}$ if $X$ is reduced and $d \geq 1$. These are the content of respectively Theorem 4 and Theorem 8 and they generalize a result in \cite{SVV} for hypersurfaces. As in that paper, we use our results to get (in the case $k = \mathbb{C}$) a lower bound on the poles of (a slight specialization of) the motivic zeta function (Theorem 9) and the topological zeta function (Theorem 10). The topological zeta function of $X$ has no pole less than $-N+(N-m+1)/2$ if $m < N$ and no pole less than $-N + (d+1)/2$ if $X$ is reduced and $d \geq 1$. We will give an example which proves that this lower bound is optimal.

In the description of our results, we started with a closed subscheme $X$ of $\mathbb{A}^N$. The results involving the number of equations can be generalized to closed subschemes $X$ of a nonsingular irreducible algebraic variety for which the ideal sheaf has locally $m$ generators. The results involving the dimension can be generalized to separated schemes of finite type. The reason for this is that we can use coverings by affine schemes.

\vspace{0,5cm}

\noindent \textbf{Acknowledgements.} The authors want to thank
Willem Veys for his support and useful remarks.

\section{Results involving the number of equations}

Let $X$ be a closed subscheme of $\mathbb{A}^N$. We denote by
$I_X$ the ideal of $k[x_1,\ldots,x_N]$ which consists of the
global sections of the ideal sheaf of $X$. We take generators
$f_1,\ldots,f_m$ of $I_X$, so $I_X=(f_1,\ldots,f_m)$. We want to
stress that $m$ is the number of polynomials which generate this
ideal.

Let $\gamma(t) \in \mathcal{L}(\mathbb{A}^N)$. We consider the
quotient $M=M(f_1,\ldots,f_m;\gamma)$ of the $k[[t]]$-module
$k[[t]]^N$ by the module generated by the rows of the matrix
\[
J =J(f_1,\ldots,f_m;\gamma) := \left( \begin{array}{ccccc}
\frac{\partial f_1}{\partial x_1}(\gamma) & \frac{\partial
f_1}{\partial x_2}(\gamma) & \cdots
& \cdots & \frac{\partial f_1}{\partial x_N}(\gamma) \\
\frac{\partial f_2}{\partial x_1}(\gamma) & \frac{\partial
f_2}{\partial x_2}(\gamma) & \cdots & \cdots & \frac{\partial
f_2}{\partial x_N}(\gamma) \\ \vdots & \vdots & & & \vdots
\\ \frac{\partial f_m}{\partial x_1}(\gamma) & \frac{\partial f_m}{\partial x_2}(\gamma) & \cdots
& \cdots & \frac{\partial f_m}{\partial x_N}(\gamma)\end{array}
\right).
\]
There exists a unique sequence $0 \leq e_1 \leq e_2 \leq \cdots
\leq e_p$ in $\mathbb{N}$, with $0\leq p \leq \min\{m,N\}$, such
that
\[
M \cong \frac{k[[t]]}{(t^{e_1})} \oplus \frac{k[[t]]}{(t^{e_2})}
\oplus \cdots \oplus \frac{k[[t]]}{(t^{e_p})} \oplus
(k[[t]])^{N-p}.
\]
We write $p=p(f_1,\ldots,f_m;\gamma)$ and
$e_i=e_i(f_1,\ldots,f_m;\gamma)$. Note that these values are
determined by the Smith normal form of $J$ and that
$t^{e_1},t^{e_2},\ldots,t^{e_p}$ are called the invariant factors
of $M$ and $J$. If $\gamma \in \mathcal{L}(X)$, one checks easily
that the module $M$ does not depend on the chosen generators of
$I_X$ and we will write $M=M(X;\gamma)$, $p=p(X;\gamma)$ and
$e_i=e_i(X;\gamma)$. Note that $p(X;\gamma)=m$ if and only if the generic point of $\gamma$ lies on the smooth part of an irreducible component of $X$ of codimension $m$, and on no other irreducible component of $X$.

For a sequence $0 \leq e_1 \leq e_2 \leq \cdots \leq e_p$ in
$\mathbb{N}$, let $\mathcal{M}^{e_1,\ldots,e_p}(f_1,\ldots,f_m)$
denote the set of all $\gamma \in \mathcal{L}(\mathbb{A}^N)$ for
which the invariant factors of $M(f_1,\ldots,f_m;\gamma)$ are
$t^{e_1},t^{e_2},\ldots,t^{e_p}$. We put
$\mathcal{L}^{e_1,\ldots,e_p}(X) := \mathcal{L}(X) \cap
\mathcal{M}^{e_1,\ldots,e_p}(f_1,\ldots,f_m)$.

Let $\theta(t) \in \mathcal{L}_n(\mathbb{A}^N)$. We consider the
quotient $M_n=M_n(f_1,\ldots,f_m;\theta)$ of the $k[[t]]$-module
$(k[[t]]/(t^{n+1}))^N$ by the module generated by the rows of the
matrix
\[
J =J(f_1,\ldots,f_m;\theta) := \left( \begin{array}{ccccc}
\frac{\partial f_1}{\partial x_1}(\theta) & \frac{\partial
f_1}{\partial x_2}(\theta) & \cdots
& \cdots & \frac{\partial f_1}{\partial x_N}(\theta) \\
\frac{\partial f_2}{\partial x_1}(\theta) & \frac{\partial
f_2}{\partial x_2}(\theta) & \cdots & \cdots & \frac{\partial
f_2}{\partial x_N}(\theta) \\ \vdots & \vdots & & & \vdots
\\ \frac{\partial f_m}{\partial x_1}(\theta) & \frac{\partial f_m}{\partial x_2}(\theta) & \cdots
& \cdots & \frac{\partial f_m}{\partial x_N}(\theta)\end{array}
\right).
\]
There exists a unique sequence $0 \leq e_1 \leq e_2 \leq \cdots
\leq e_N \leq n+1$ in $\mathbb{N}$ such that
\[
M_n \cong \frac{k[[t]]}{(t^{e_1})} \oplus \frac{k[[t]]}{(t^{e_2})}
\oplus \cdots \oplus \frac{k[[t]]}{(t^{e_N})}.
\]
We write $e_i=e_i(f_1,\ldots,f_m;\theta)$. If $\theta \in
\mathcal{L}_n(X)$, then the module $M_n$ does not depend on the
chosen generators of $I_X$ and we will write $M_n=M_n(X;\theta)$
and $e_i=e_i(X;\theta)$. Let $g$ be $n/2$ if $n$ is even and
$(n+1)/2$ if $n$ is odd. For $\theta \in \mathcal{L}_n(X)$, the
largest $i \in \{0,\ldots,N\}$ for which $e_i(X;\theta) < g$ is
denoted by $b(X;\theta)$. Here, we have put $e_0(X;\theta)=0$.

For a sequence $0 \leq e_1 \leq e_2 \leq \cdots \leq e_b<g$ in
$\mathbb{N}$ with $b \leq N$, let
$\mathcal{M}_n^{e_1,\ldots,e_b}(f_1,\ldots,f_m)$ denote the set of
all $\theta \in \mathcal{L}_n(\mathbb{A}^N)$ satisfying
$e_i(f_1,\ldots,f_m;\theta)=e_i$ for $i=1,\ldots,b$ and
$e_{b+1}(f_1,\ldots,f_m;\theta) \geq g$ if $b < N$. We put
$\mathcal{L}_n^{e_1,\ldots,e_b}(X) := \mathcal{L}_n(X) \cap
\mathcal{M}_n^{e_1,\ldots,e_b}(f_1,\ldots,f_m)$.

\vspace{0,5cm}

\noindent \textbf{Proposition 1.} Suppose that $m \leq N$ and let
$0 \leq e_1 \leq e_2 \leq \cdots \leq e_m$ be an increasing
sequence of integers.

(a) For $l \geq e_m$ and $i \geq e_m$, we have
\[
\pi_l^{l+i}(\mathcal{L}_{l+i}^{e_1,\ldots,e_m}(X)) =
\pi_l(\mathcal{L}^{e_1,\ldots,e_m}(X)),
\]
and this set is locally closed.

(b) For $\gamma \in \pi_{e_m}(\mathcal{L}^{e_1,\ldots,e_m}(X))$,
there exists an open subset $U_{e_m}$ of
$\pi_{e_m}(\mathcal{L}^{e_1,\ldots,e_m}(X))$ containing $\gamma$
such that if we denote by $U_i$ the inverse image of $U_{e_m}$
under the truncation map
\[
\pi_i(\mathcal{L}^{e_1,\ldots,e_m}(X)) \rightarrow
\pi_{e_m}(\mathcal{L}^{e_1,\ldots,e_m}(X))
\]
for $i >e_m$, we have for $l \geq e_m$ that the truncation map
$U_{l+1} \rightarrow U_l$ is a trivial fibration with fibre
$\mathbb{A}^{N-m}$.

(c) For $n >l \geq e_m$, the truncation map
\[
\pi_n(\mathcal{L}^{e_1,\ldots,e_m}(X)) \rightarrow
\pi_l(\mathcal{L}^{e_1,\ldots,e_m}(X))
\]
is a locally trivial fibration with fibre
$\mathbb{A}^{(N-m)(n-l)}$.

\vspace{0,5cm}

\noindent \textbf{Proof.} Let $J$ be the $(m \times N)$-matrix
with component in the $i$-th row and $j$-th column equal to
\[
\frac{\partial f_i}{\partial
x_j}(a_{1,0}+a_{1,1}t+\cdots+a_{1,e_m}t^{e_m}, \ldots,
a_{N,0}+a_{N,1}t+\cdots+a_{N,e_m}t^{e_m}),
\]
which is an element of
$\mathcal{O}(\mathcal{L}_{e_m}(\mathbb{A}^N))[[t]]$ and where
$(a_{i,j})_{1 \leq i \leq N;0 \leq j \leq e_m}$ are the
coordinates on $\mathcal{L}_{e_m}(\mathbb{A}^N)$. By using row
operations, one checks that there exists
\begin{enumerate} \vspace{-2mm}
\item an open neighborhood $V$ of $\gamma$ in
$\mathcal{L}_{e_m}(\mathbb{A}^N)$ and \vspace{-2mm}
\item an $(m \times m)$-matrix $A$ with components in $\mathcal{O}(V)[[t]]$,
which is invertible in the ring of $m \times m$ matrices over
$\mathcal{O}(V)[[t]]$
\end{enumerate}
such that
\begin{enumerate} \vspace{-2mm}
\item $V \cap
\pi_{e_m}(\mathcal{M}^{e_1,\ldots,e_m}(f_1,\ldots,f_m))$ is a
closed subset of $V$ and \vspace{-2mm}
\item on $V \cap \pi_{e_m}(\mathcal{M}^{e_1,\ldots,e_m}(f_1,\ldots,f_m))$
and after renumbering the variables
\begin{eqnarray} \label{matrixeq}
A \cdot J = \left(
\begin{array}{ccccccc} t^{e_1} & h_{1,2} & h_{1,3} & \cdots & h_{1,m} & \cdots & h_{1,N} \\ 0 &
t^{e_2} & h_{2,3} & \cdots & h_{2,m} & \cdots & h_{2,N} \\ 0 & 0 &
t^{e_3} & \cdots & h_{3,m} & \cdots & h_{3,N} \\ \vdots & \vdots &
& \ddots & \vdots & \vdots & \vdots \\ 0 & 0&  \cdots & 0 &
t^{e_m} & \cdots & h_{m,N}
\end{array} \right),
\end{eqnarray}
\end{enumerate}
where there are zeroes under the diagonal and where every
$h_{i,j}$ is a formal power series in $t$ of order at least $e_i$
with coefficients restrictions of regular functions on $V$ to $V
\cap \pi_{e_m}(\mathcal{M}^{e_1,\ldots,e_m}(f_1,\ldots,f_m))$.

We have to consider the system of equations
\[ \left\{ \begin{array}{c}
f_1(a_{1,0}+a_{1,1}t+a_{1,2}t^2+\cdots,\ldots,a_{N,0}+a_{N,1}t+a_{N,2}t^2+\cdots)
= 0 \\
f_2(a_{1,0}+a_{1,1}t+a_{1,2}t^2+\cdots,\ldots,a_{N,0}+a_{N,1}t+a_{N,2}t^2+\cdots)
= 0 \\ \vdots \\
f_m(a_{1,0}+a_{1,1}t+a_{1,2}t^2+\cdots,\ldots,a_{N,0}+a_{N,1}t+a_{N,2}t^2+\cdots)
= 0
\end{array} \right. \]
Let now $l$ be an integer which is at least $e_m$. We denote
\begin{eqnarray*}
a & = &
(a_{1,0}+a_{1,1}t+\cdots+a_{1,l}t^l,\ldots,a_{N,0}+a_{N,1}t+\cdots+a_{N,l}t^l)
\\ z & = & (z_1,\ldots,z_N) \\ & = & (a_{1,l+1}+a_{1,l+2}t+a_{1,l+3}t^2+\cdots,\ldots,
a_{N,l+1}+a_{N,l+2}t+a_{N,l+3}t^2+\cdots)
\end{eqnarray*}
The previous system of equations is the same as
\[ \left\{ \begin{array}{c}
f_1(a) + \frac{\partial f_1}{\partial x_1}(a)t^{l+1}z_1 + \cdots +
\frac{\partial f_1}{\partial x_N}(a)t^{l+1}z_N + t^{2l+2}(\cdots)
= 0  \\ f_2(a) + \frac{\partial f_2}{\partial x_1}(a)t^{l+1}z_1 +
\cdots + \frac{\partial
f_2}{\partial x_N}(a)t^{l+1}z_N + t^{2l+2}(\cdots) = 0 \\
\vdots \\ f_m(a) + \frac{\partial f_m}{\partial x_1}(a)t^{l+1}z_1
+ \cdots + \frac{\partial f_m}{\partial x_N}(a)t^{l+1}z_N +
t^{2l+2}(\cdots) = 0
\end{array} \right. \]
Let $J_l$ be the $(m \times N)$-matrix with component in the
$i$-th row and $j$-th column equal to $(\partial f_i/\partial
x_j)(a)$. We now multiply our system of equations with $A$, and
get on $(\pi_{e_m})^{-1}(V)$ the equivalent system of equations
\begin{eqnarray} \label{eqsystem}
\left( \begin{array}{c} h_1 \\ h_2 \\ \vdots \\ h_m \end{array}
\right) + A \cdot J_l \left( \begin{array}{c} t^{l+1}z_1 \\ t^{l+1}z_2 \\
\vdots \\ t^{l+1}z_N \end{array} \right) + \left(
\begin{array}{c} t^{2l+2}u_1 \\ t^{2l+2}u_2 \\ \vdots \\ t^{2l+2}u_m \end{array}
\right) = \left( \begin{array}{c} 0 \\ 0 \\ \vdots \\ 0
\end{array} \right),
\end{eqnarray}
where $h_i \in \mathcal{O}((\pi^l_{e_m})^{-1}(V))[[t]]$ and where
$u_i \in
(\mathcal{O}((\pi^l_{e_m})^{-1}(V))[[t]])[z_1,\ldots,z_N]$ has
multiplicity at least 2 as polynomial in the $z_1,\ldots,z_N$. We
use now that $A \cdot J_l \equiv A \cdot J \mbox{ mod } t^{e_m+1}$
and that the equality of matrices in (\ref{matrixeq}) holds on
$(\pi_{e_m})^{-1}(V) \cap
\mathcal{M}^{e_1,\ldots,e_m}(f_1,\ldots,f_m)$, and obtain
necessary conditions for an element of $(\pi_{e_m})^{-1}(V) \cap
\mathcal{M}^{e_1,\ldots,e_m}(f_1,\ldots,f_m)$ to be in
$\mathcal{L}^{e_1,\ldots,e_m}(X)$, namely
\begin{eqnarray} \label{conditions}
\left\{ \begin{array}{c} h_1 \equiv 0 \mbox{ mod } t^{l+e_1+1}
\\ h_2 \equiv 0 \mbox{ mod } t^{l+e_2+1} \\ \vdots \\ h_m \equiv 0 \mbox{ mod } t^{l+e_m+1} \end{array}
\right..
\end{eqnarray}
Let us denote the subset of $(\pi^l_{e_m})^{-1}(V) \cap
\pi_l(\mathcal{M}^{e_1,\ldots,e_m}(f_1,\ldots,f_m))$ given by
these conditions by
$\mathcal{N}_l^{e_1,\ldots,e_m}(f_1,\ldots,f_m)$. By using (\ref{eqsystem}), one can check that the truncation map
$\mathcal{N}_{l+1}^{e_1,\ldots,e_m}(f_1,\ldots,f_m) \rightarrow
\mathcal{N}_l^{e_1,\ldots,e_m}(f_1,\ldots,f_m)$ is a trivial
fibration with fibre $\mathbb{A}^{N-m}$. We can continue in this
way, so that every $l$-jet in
$\mathcal{N}_l^{e_1,\ldots,e_m}(f_1,\ldots,f_m)$ is liftable to an
arc in $\mathcal{L}^{e_1,\ldots,e_m}(X)$. Consequently,
\[
(\pi^l_{e_m})^{-1}(V) \cap \pi_l(\mathcal{L}^{e_1,\ldots,e_m}(X))
= \mathcal{N}_l^{e_1,\ldots,e_m}(f_1,\ldots,f_m).
\]
The statements in (a) follow now from the fact that the $l$-jets
in the closed subset $(\pi^l_{e_m})^{-1}(V) \cap
\pi_l(\mathcal{M}^{e_1,\ldots,e_m}(f_1,\ldots,f_m))$ of
$(\pi^l_{e_m})^{-1}(V)$ which can be lifted to an arc in
$\mathcal{L}^{e_1,\ldots,e_m}(X)$ are given by the conditions in
(\ref{conditions}). We obtain (b) if we put $U_l :=
\mathcal{N}_l^{e_1,\ldots,e_m}(f_1,\ldots,f_m)$. Part (c) follows
in a straightforward way from part (b). $\qquad \Box$

\vspace{0,5cm}

In proposition 2 and 3, we will look at the space of $n$-jets. As
before, $g$ will be $n/2$ if $n$ is even and $(n+1)/2$ if $n$ is
odd.

\vspace{0,5cm}

\noindent \textbf{Proposition 2.} Suppose that $m \leq N$ and let
$0 \leq e_1 \leq e_2 \leq \cdots \leq e_m < g$ be an increasing
sequence of integers. Put $e=e_1+\cdots+e_m$.

(a) For $l \in \{e_m,e_m+1,\ldots,n\}$, we have that
\[
\pi^n_l(\mathcal{L}_n^{e_1,\ldots,e_m}(X))
\]
is locally closed.

(b) For $\gamma \in
\pi^n_{e_m}(\mathcal{L}_n^{e_1,\ldots,e_m}(X))$, there exists an
open subset $U_{e_m}$ of
$\pi^n_{e_m}(\mathcal{L}_n^{e_1,\ldots,e_m}(X))$ containing
$\gamma$ such that if we denote by $U_i$ the inverse image of
$U_{e_m}$ under the truncation map
\[
\pi^n_i(\mathcal{L}_n^{e_1,\ldots,e_m}(X)) \rightarrow
\pi^n_{e_m}(\mathcal{L}_n^{e_1,\ldots,e_m}(X))
\]
for $i \in \{e_m+1,\ldots,n\}$, we have for $l \in
\{e_m,\ldots,n-1\}$ that the truncation map $U_{l+1} \rightarrow
U_l$ is a trivial fibration with fibre $\mathbb{A}^{N-a}$, where
$a$ is the largest number of $\{0,\ldots,m\}$ satisfying
$e_a<n-l$.

(c) For $n>l\geq e_m$, the truncation map
\[
\mathcal{L}_n^{e_1,\ldots,e_m}(X) \rightarrow
\pi^n_l(\mathcal{L}_n^{e_1,\ldots,e_m}(X))
\]
is a locally trivial fibration with fibre an affine space. If $l
\leq n-e_m$, the dimension of this fibre is equal to
$(N-m)(n-l)+e$. If $l=e_m$ and $m < N$, this dimension is at least
$\ulcorner (N-m+1)(n/2) \urcorner$.

\vspace{0,5cm}

\noindent \textbf{Proof.} Part (a) and (b) can be proved in an
analogous way as in the previous proposition. Instead of working
with power series in $t$, we have to work with power series
truncated modulo $t^{n+1}$. We consider now part (c). It follows
from (b) that the map
\[
\mathcal{L}_n^{e_1,\ldots,e_m}(X) \rightarrow
\pi^n_l(\mathcal{L}_n^{e_1,\ldots,e_m}(X))
\]
is a locally trivial fibration with fibre an affine space. If $l
\leq n-e_m$, the dimension of this fibre is equal to
\begin{eqnarray*}
\lefteqn{Ne_1+(N-1)(e_2-e_1)+ \cdots +
(N-m+1)(e_m-e_{m-1})+(N-m)(n-e_m-l) } \\ & = &
e_1+e_2+\cdots+e_m+(N-m)e_m+(N-m)(n-e_m-l) \hspace{3cm} \\
& = & e+(N-m)(n-l).
\end{eqnarray*}
Our lower bound when $l=e_m$ and $m<N$ follows now from the
calculation
\begin{eqnarray*}
\lefteqn{e+(N-m)(n-e_m)} \\ & = & e_1+e_2+\cdots+e_{m-1}+2e_m-n+(N-m+1)(n-e_m) \\
& = & e_1+e_2+\cdots+e_{m-1}+2e_m-n+(N-m+1)(n/2-e_m)+(N-m+1)(n/2) \\
& \geq & e_1+e_2+\cdots+e_{m-1}+2e_m-n+n-2e_m+(N-m+1)(n/2) \\ & =
& e_1+e_2+\cdots+e_{m-1}+(N-m+1)(n/2) \\ & \geq & (N-m+1)(n/2)
\end{eqnarray*}
and from the fact that the dimension of the fibre is an integer.
$\qquad \Box$

\vspace{0,5cm}

\noindent \textbf{Proposition 3.} Let $0 \leq e_1 \leq e_2 \leq
\cdots \leq e_b < g$ be an increasing sequence of integers.
Suppose that $b \leq m$ and $b \leq N$.

(a) For $l \in \{n-g,\ldots,n\}$, we have that
\[
\pi^n_l(\mathcal{L}^{e_1,\ldots,e_b}_n(X))
\]
is locally closed.

(b) For $\gamma \in
\pi^n_{n-g}(\mathcal{L}_n^{e_1,\ldots,e_b}(X))$, there exists an
open subset $U_{n-g}$ of
$\pi^n_{n-g}(\mathcal{L}_n^{e_1,\ldots,e_b}(X))$ containing
$\gamma$ such that if we denote by $U_i$ the inverse image of
$U_{n-g}$ under the truncation map
\[
\pi^n_i(\mathcal{L}_n^{e_1,\ldots,e_b}(X)) \rightarrow
\pi^n_{n-g}(\mathcal{L}_n^{e_1,\ldots,e_b}(X))
\]
for $i \in \{n-g+1,\ldots,n\}$, we have for $l \in
\{n-g,\ldots,n-1\}$ that the truncation map $U_{l+1} \rightarrow
U_l$ is a trivial fibration with fibre $\mathbb{A}^{N-a}$, where
$a$ is the largest number of $\{0,\ldots,b\}$ satisfying
$e_a<n-l$.

(c) For $n > l \geq n-g$, the truncation map
\[
\mathcal{L}_n^{e_1,\ldots,e_b}(X) \rightarrow
\pi^n_l(\mathcal{L}_n^{e_1,\ldots,e_b}(X))
\]
is a locally trivial fibration with fibre an affine space. If
$l=n-g$, the dimension of this fibre is at least $\ulcorner
(N-b)(n/2) \urcorner$ and if moreover $e_1+\cdots+e_b \geq n/2$,
then it is at least $\ulcorner (N-b+1)(n/2) \urcorner$.

(d) Let $\theta \in \mathcal{L}_n^{e_1,\ldots,e_b}(X)$. Let $M$ be
a closed subscheme of $\mathbb{A}^N$ defined by $b$ elements of
$I_X$ such that $\theta \in \mathcal{L}_n^{e_1,\ldots,e_b}(M)$. If
$\alpha \in \mathcal{L}_n^{e_1,\ldots,e_b}(M)$ satisfies
$\pi^n_{n-g}(\alpha)=\pi^n_{n-g}(\theta)$, then $\alpha \in
\mathcal{L}_n^{e_1,\ldots,e_b}(X)$.

\vspace{0,5cm}

\noindent \textbf{Proof.} Part (a),(b) and (d) are again proved
analogously as before. We adapt the first part of the proof of
Proposition 1, the remaining part of the adaptation  is left to
the reader. Let $J$ be the $(m \times N)$-matrix with component in
the $i$-th row and $j$-th column equal to
\[
\frac{\partial f_i}{\partial
x_j}(a_{1,0}+a_{1,1}t+\cdots+a_{1,n-g}t^{n-g}, \ldots,
a_{N,0}+a_{N,1}t+\cdots+a_{N,n-g}t^{n-g})
\]
considered as an element of
$\mathcal{O}(\mathcal{L}_{n-g}(\mathbb{A}^N))[[t]]/(t^{n+1})$.
There exists
\begin{enumerate} \vspace{-2mm}
\item an open neighborhood $V$ of $\gamma$ in
$\mathcal{L}_{n-g}(\mathbb{A}^N)$ and \vspace{-2mm}
\item an $(m \times m)$-matrix $A$ with components in $\mathcal{O}(V)[[t]]/(t^{n+1})$,
which is invertible in the ring of $m \times m$ matrices over
$\mathcal{O}(V)[[t]]/(t^{n+1})$
\end{enumerate}
such that
\begin{enumerate} \vspace{-2mm}
\item $V \cap
\pi^n_{n-g}(\mathcal{M}_n^{e_1,\ldots,e_b}(f_1,\ldots,f_m))$ is a closed subset of $V$ and \vspace{-2mm}
\item on $V \cap \pi^n_{n-g}(\mathcal{M}_n^{e_1,\ldots,e_b}(f_1,\ldots,f_m))$
and after renumbering the variables
\begin{eqnarray*}
A \cdot J = \left(
\begin{array}{cccccccccc} t^{e_1} & h_{1,2} & h_{1,3} & \cdots & h_{1,b-1} & h_{1,b} & h_{1,b+1} & \cdots & h_{1,N} \\ 0 &
t^{e_2} & h_{2,3} & \cdots & h_{2,b-1} & h_{2,b} & h_{2,b+1} &
\cdots & h_{2,N} \\ 0 & 0 & t^{e_3} & \cdots & h_{3,b-1} & h_{3,b} & h_{3,b+1} & \cdots & h_{3,N} \\
\vdots & \vdots & & \ddots & & \vdots & \vdots & & \vdots
\\ 0 & 0 & & & t^{e_{b-1}} & h_{b-1,b} & h_{b-1,b+1} & \cdots & h_{b-1,N} \\
0 & 0& \cdots & \cdots & 0 & t^{e_b} & h_{b,b+1} & \cdots & h_{b,N} \\
0 & 0& \cdots & \cdots & 0 & 0& h_{b+1,b+1} & \cdots & h_{b+1,N} \\
\vdots & \vdots & & & \vdots & \vdots & \vdots & & \vdots \\ 0 &
0& \cdots & \cdots & 0 & 0 & h_{m,b+1} & \cdots & h_{m,N}
\end{array} \right),
\end{eqnarray*}
\end{enumerate}
where every $h_{i,j}$ is obtained by restricting the coefficients
of an element of $\mathcal{O}(V)[[t]]/(t^{n+1})$ to $V \cap
\pi^n_{n-g}(\mathcal{M}_n^{e_1,\ldots,e_b}(f_1,\ldots,f_m))$, and
where $h_{i,j}$ has order at least $e_i$ if $b,j \geq i$ and at
least $g$ if $i,j \geq b+1$.

We finally check the lower bound for the dimension of the fibre in
part (c) in the case $l=n-g$. It follows from part (b) that the
dimension of the fibre is equal to
\begin{eqnarray*}
\lefteqn{Ne_1+(N-1)(e_2-e_1)+ \cdots +
(N-b+1)(e_b-e_{b-1})+(N-b)(g-e_b)}  \\ & = &
e_1+e_2+\cdots+e_b+(N-b)g. \hspace{7cm}
\end{eqnarray*}
This is at least $\ulcorner (N-b)(n/2) \urcorner$. If moreover
$e_1+\cdots+e_b \geq n/2$, it is at least $\ulcorner (N-b+1)(n/2)
\urcorner$. $\qquad \Box$

\vspace{0,5cm}

\noindent \textbf{Theorem 4.} Let $X$ be a closed subscheme of
$\mathbb{A}^N$ and suppose that $I_X$ is generated by $m$
polynomials, with $m<N$. Then $[\mathcal{L}_n(X)]$ is a multiple
of $\mathbb{L}^{\ulcorner (N-m+1)(n/2) \urcorner}$ in
$K_0(\mbox{Var}_{k})$ for all $n \in \mathbb{Z}_{\geq 0}$.

\vspace{0,5cm}

\noindent \textbf{Proof.} This follows from Proposition 2(c) and
Proposition 3(c). $\qquad \Box$

\section{Results involving the dimension}

Let $X$ be a closed subscheme of $\mathbb{A}^N$. Let $d$ be the
minimum of the dimensions of the irreducible components of $X$ and
put $r:=N-d$. We denote by $\mbox{Jac}_X$ the scheme defined by
the ideal of $k[x_1,\ldots,x_N]$ generated by $I_X$ and all $r$ by
$r$ minors of matrices of the form $\partial(g_1,\ldots,g_r) /
\partial(x_1,\ldots,x_N)$, where $g_1,\ldots,g_r$ are polynomials in
$I_X$. Because all $r+1$ by $r+1$ minors of matrices of the form
$\partial(g_1,\ldots,g_{r+1}) / \partial(x_1,\ldots,x_N)$, with
$g_1,\ldots,g_{r+1} \in I_X$, lie in the radical of $I_X$, we get that
$p(X;\gamma) \leq r$ for $\gamma \in \mathcal{L}(X)$. Moreover,
for $\gamma \in \mathcal{L}(X)$, we get $\gamma \not\in
\mathcal{L}(\mbox{Jac}_X)$ if and only if $p(X;\gamma)=r$, and if
this is the case, $\gamma$ lies on an irreducible component of $X$
of codimension $r$. Note that $p(X;\gamma)$ is defined in the beginning of Section 2, and that it is the number of elements $e_i(X;\gamma)$. For $\gamma \in \mathcal{L}(X) \setminus
\mathcal{L}(\mbox{Jac}_X)$, we put $e(X;\gamma) =
\mbox{ord}_{\gamma}(\mbox{Jac}_X) = \sum_{i=1}^r e_i(X;\gamma)$.
For $\theta \in \mathcal{L}_n(X)$, we put $e(X;\theta):=
\sum_{i=1}^r e_i(X;\theta)$.

Suppose that $I_X=(f_1,\ldots,f_m)$. Fix polynomials
$g_1,\ldots,g_r$ in $I_X$ and let $M$ be the scheme defined by the
ideal $(g_1,\ldots,g_r)$. Suppose that the minimum of the
dimensions of the irreducible components of $M$ is also $d$. Let
$X'$ be the reduced closed subscheme of $\mathbb{A}^N$ whose
irreducible components are the irreducible components of $M$ which are
not a component of $X$ and let $I_{X'}=(f_1',\ldots,f_s')$. Note
that set-theoretically $M=X \cup X'$, which implies
$\mathcal{L}(M)=\mathcal{L}(X) \cup \mathcal{L}(X')$. We endow $X
\cap X'$ with the reduced induced subscheme structure. Let $I_{X
\cap X'} = (h_1,\ldots,h_t)$. By Hilbert's Nullstellensatz, there
exists $c \in \mathbb{N}$ such that $h_i^c \in
(f_1,\ldots,f_m,f_1',\ldots,f_s')$ for all $i \in \{1,\ldots,t\}$. Note that $c$ depends on $M$, and therefore we will write $c=c(M)$ if necessary.

\vspace{0,5cm}

\noindent \textbf{Lemma A.} Let $n \in \mathbb{N}$. If $\theta \in
\mathcal{L}_{cn}(X) \cap \mathcal{L}_{cn}(X')$, then we have
$\pi^{cn}_n(\theta) \in \mathcal{L}_{n}(X \cap X')$ and thus
$\mbox{ord}_{\theta}(\mbox{Jac}_M)>n$.

\vspace{0,5cm}

\noindent \textbf{Proof.} Let $\theta \in \mathcal{L}_{cn}(X) \cap
\mathcal{L}_{cn}(X')$. For every $i$, we have that $h_i^c(\theta)
\equiv 0 \mbox{ mod }t^{cn+1}$ and thus $h_i(\theta) \equiv 0
\mbox{ mod } t^{n+1}$. Consequently, $\pi^{cn}_n(\theta) \in
\mathcal{L}_{n}(X \cap X')$. Because
$I_{\mbox{\footnotesize{Jac}}_M} \subset I_{X \cap X'}$, we obtain
$\mbox{ord}_{\theta}(\mbox{Jac}_M)>n$. $\qquad \Box$

\vspace{0,5cm}

\noindent \textbf{Lemma B.} Let $\rho \in \mathcal{L}(M) \setminus
\mathcal{L}(\mbox{Jac}_M)$ and put $e=e(M;\rho)$. If $\pi_{ce}(\rho)
\in \mathcal{L}_{ce}(X)$, then $\rho \in \mathcal{L}(X)$.

\vspace{0,5cm}

\noindent \textbf{Proof.} Since $\mbox{ord}_{\rho}(\mbox{Jac}_M)=e$,
we have by Lemma A that $\pi_{ce}(\rho) \not\in \mathcal{L}_{ce}(X)
\cap \mathcal{L}_{ce}(X')$. Since $\pi_{ce}(\rho) \in
\mathcal{L}_{ce}(X)$, we obtain that $\pi_{ce}(\rho) \not\in
\mathcal{L}_{ce}(X')$ and thus $\rho \not\in \mathcal{L}(X')$.
Because $\rho \in \mathcal{L}(M) = \mathcal{L}(X) \cup
\mathcal{L}(X')$, this implies $\rho \in \mathcal{L}(X)$. $\qquad
\Box$

\vspace{0,5cm}

\noindent \textbf{Proposition 5.} (a) Let $\gamma \in
\mathcal{L}(X) \setminus \mathcal{L}(\mbox{Jac}_X)$. For $l \in
\mathbb{N}$, we define
\[
A_l := \{ \rho \in\mathcal{L}(\mathbb{A}^N) \mid \pi_l(\rho) =
\pi_l(\gamma) \}.
\]
Put $e_i = e_i(X;\gamma)$ for $i=1,\ldots,p(X;\gamma)=r$. Take a
closed subscheme $M$ defined by $r$ elements $g_1,\ldots,g_r$ of
$I_X$ such that $\gamma \in \mathcal{L}^{e_1,\ldots,e_r}(M)$.
Then,
\begin{enumerate} \vspace{-3mm}
\item[(i)] for $l \geq e_r$, we have $\mathcal{L}(M) \cap A_l = \mathcal{L}(X) \cap
A_l$ and

\vspace{-3mm}
\item[(ii)] for $n > l \geq e_r$, we have
$\pi_n(\mathcal{L}(X) \cap A_l) \cong \mathbb{A}^{(n-l)d}$.
\end{enumerate}
(b) For $l \geq e_r$, the set
$\pi_l(\mathcal{L}^{e_1,\ldots,e_r}(X))$ is constructible and for
$n
>l \geq e_r$, the truncation map
\[
\pi_n(\mathcal{L}^{e_1,\ldots,e_r}(X)) \rightarrow
\pi_l(\mathcal{L}^{e_1,\ldots,e_r}(X))
\]
is a piecewise trivial fibration with fibre $\mathbb{A}^{(n-l)d}$.

\vspace{0,5cm}

\noindent \textbf{Proof.} We explained in the beginning of this
section that $p(X;\gamma)=r$. Because $\gamma \in
\mathcal{L}^{e_1,\ldots,e_r}(M)$, we get that $\gamma$ lies on an
irreducible component $M_{\gamma}$ of $M$ of codimension at least
$r$. Because $X \subset M$, the codimension of $M_{\gamma}$ has to
be equal to $r$, and every other irreducible component of $M$ has
codimension at most $r$. Note also that $M_{\gamma}$ has to be an
irreducible component of $X$.

(a.i) We only have to prove that $\mathcal{L}(M) \cap A_{e_r}
\subset \mathcal{L}(X)$. Note that for every $\rho \in
\mathcal{L}(M) \cap A_{e_r}$, we have that
$\pi_{e_r}(\gamma)=\pi_{e_r}(\rho)$ and $e(M;\rho) = e(M;\gamma)
=:e$.

It follows from Lemma A that the union
\[
\pi_{ce}(\mathcal{L}(M) \cap A_{e_r}) = \left( \mathcal{L}_{ce}(X)
\cap \pi_{ce}(\mathcal{L}(M) \cap A_{e_r}) \right) \cup \left(
\mathcal{L}_{ce}(X') \cap \pi_{ce}(\mathcal{L}(M) \cap A_{e_r})
\right)
\]
is actually a disjoint union. It follows from Proposition 1(c)
that $\pi_{ce}(\mathcal{L}(M) \cap A_{e_r})$ is an irreducible
closed subset of $\mathcal{L}_{ce}(\mathbb{A}^N)$ isomorphic to
$\mathbb{A}^{(ce-e_r)d}$. We also have that $\mathcal{L}_{ce}(X)
\cap \pi_{ce}(\mathcal{L}(M) \cap A_{e_r})$ and
$\mathcal{L}_{ce}(X') \cap \pi_{ce}(\mathcal{L}(M) \cap A_{e_r})$
are closed subsets of $\mathcal{L}_{ce}(\mathbb{A}^N)$. Because
$\pi_{ce}(\gamma) \in \mathcal{L}_{ce}(X) \cap
\pi_{ce}(\mathcal{L}(M) \cap A_{e_r})$, we obtain
\[
\pi_{ce}(\mathcal{L}(M) \cap A_{e_r}) = \mathcal{L}_{ce}(X) \cap
\pi_{ce}(\mathcal{L}(M) \cap A_{e_r})
\]
and of course also $\mathcal{L}_{ce}(X') \cap
\pi_{ce}(\mathcal{L}(M) \cap A_{e_r}) = \emptyset$.

Take $\rho \in \mathcal{L}(M) \cap A_{e_r}$. We have proven that
$\pi_{ce}(\rho) \in \mathcal{L}_{ce}(X)$. Hence, $\rho \in
\mathcal{L}(X)$ by Lemma B.

(a.ii) This follows from the previous part and from Proposition
1(c).

(b) Put $e=e_1+\cdots+e_r$. Suppose that $I_X = (f_1,\ldots,f_m)$. 
There is a closed subscheme $M_i$ associated to every $r$ elements of $\{f_1,\ldots,f_m\}$, and we consider only those $M_i$ which have an irreducible component of dimension $d$. Clearly, there are only a finite number of such subschemes $M_i$. Put $c = \max \{ c(M_i) \}$.

We prove now that $\pi_l(\mathcal{L}^{e_1,\ldots,e_r}(X)) =
\pi^{ce+l}_l(\mathcal{L}^{e_1,\ldots,e_r}_{ce+l}(X))$, which
implies that $\pi_l(\mathcal{L}^{e_1,\ldots,e_r}(X))$ is
constructible. Let $\theta \in
\mathcal{L}^{e_1,\ldots,e_r}_{ce+l}(X)$. Take a closed subscheme
$M_i$ of $\mathbb{A}^N$ defined by $r$ elements of $\{f_1,\ldots,f_m\}$ such that
$\theta \in \mathcal{L}^{e_1,\ldots,e_r}_{ce+l}(M_i)$. By
Proposition (1.a), we can lift $\pi^{ce+l}_{ce+l-e_r}(\theta)$ to
an arc $\rho \in \mathcal{L}^{e_1,\ldots,e_r}(M_i)$. Because
$\pi_{ce+l-e_r}(\rho) = \pi^{ce+l}_{ce+l-e_r}(\theta) \in
\mathcal{L}_{ce+l-e_r}(X)$, we get by Lemma B
that $\rho \in \mathcal{L}(X)$. Our statement follows now from
$\pi_l(\rho) = \pi^{ce+l}_l(\theta)$ and $\rho \in \mathcal{L}^{e_1,\ldots,e_r}(X)$ .

Finally, we proof the statement about the truncation map. Because every $\gamma \in \mathcal{L}^{e_1,\ldots,e_r}(X)$ is contained in at
least one $\mathcal{L}^{e_1,\ldots,e_r}(M_i)$, we obtain
\[
\mathcal{L}^{e_1,\ldots,e_r}(X) \subset \bigcup_i
\mathcal{L}^{e_1,\ldots,e_r}(M_i).
\]
Consequently, we have a cover
\[
\mathcal{L}^{e_1,\ldots,e_r}(X) = \bigcup_i \left(
\mathcal{L}^{e_1,\ldots,e_r}(M_i) \cap
\mathcal{L}^{e_1,\ldots,e_r}(X) \right)
\]
of $\mathcal{L}^{e_1,\ldots,e_r}(X)$. The cover
\[
\pi_l(\mathcal{L}^{e_1,\ldots,e_r}(X)) = \bigcup_i \pi_l\left(
\mathcal{L}^{e_1,\ldots,e_r}(M_i) \cap
\mathcal{L}^{e_1,\ldots,e_r}(X) \right)
\]
is an open cover of $\pi_l(\mathcal{L}^{e_1,\ldots,e_r}(X))$ by a
reasoning in the spirit of what we said in the beginning of the
proof of Proposition 1. If we apply Proposition 1(c) to $M_i$ and
use part (a.i) of this proposition, we get that the truncation map
\[
\pi_n(\mathcal{L}^{e_1,\ldots,e_r}(M_i) \cap
\mathcal{L}^{e_1,\ldots,e_r}(X)) \rightarrow
\pi_l(\mathcal{L}^{e_1,\ldots,e_r}(M_i) \cap
\mathcal{L}^{e_1,\ldots,e_r}(X))
\]
is a piecewise trivial fibration with fibre an affine space of
dimension $(n-l)d$. The result follows. $\qquad \Box$

\vspace{0,5cm}

\noindent \textbf{Proposition 6.} (a) For $\theta \in
\mathcal{L}_n(X)$, we have $b(X;\theta) \leq r$. \\ (b) For $n>l
\geq n-g$, the set $\pi^n_l(\mathcal{L}_n^{e_1,\ldots,e_b}(X))$ is
constructible and the truncation map
\[
\mathcal{L}_n^{e_1,\ldots,e_b}(X) \rightarrow
\pi^n_l(\mathcal{L}_n^{e_1,\ldots,e_b}(X))
\]
is a piecewise trivial fibration with fibre an affine space. If
$l=n-g$, the dimension of this fibre is at least $\ulcorner dn/2
\urcorner$. If moreover $b<r$ or $b=r$ and $e_1+\cdots+e_r \geq
n/2$, then this dimension is at least $\ulcorner (d+1)(n/2)
\urcorner$.

\vspace{0,5cm}

\noindent \textbf{Proof.} (a) Put $e_i=e_i(X;\theta)$ for $i \in
\{1,\ldots,N\}$. Suppose that $e_{r+1} < n/2$. Take
$g_1,\ldots,g_{r+1} \in I_X$ such that
$e_i(g_1,\ldots,g_{r+1};\theta)=e_i$ for $i=1,\ldots,r+1$. Let $Y$
be the scheme defined by the ideal $(g_1,\ldots,g_{r+1})$.
According to Proposition 1(a), $\pi^n_{n-e_{r+1}}(\theta)$ can be
lifted to an arc $\gamma \in \mathcal{L}^{e_1,\ldots,e_{r+1}}(Y)$.
But then, $\gamma$ has to lie on an irreducible component of $Y$
of codimension at least $r+1$ and this is in contradiction with
the definition of $r$ because $X \subset Y$.

(b) This follows straightforward from Proposition 3(c) and the
previous part of this proposition. $\qquad \Box$

\vspace{0,5cm}

In the remaining part of this section, we suppose that $X$ is
reduced and of pure dimension $d$. Suppose that
$I_X=(f_1,\ldots,f_m)$. Take an invertible matrix $(a_{i,j}) \in
GL_m(k)$ and denote $F_j = \sum_{i=1}^m a_{i,j}f_i$ for $1 \leq j
\leq m$. Then, $I_X=(F_1,\ldots,F_m)$. Assume that the closed
subscheme $M$ defined by $I_M=(F_1,\ldots,F_r)$ has the properties
\begin{enumerate} \vspace{-2mm}
\item all irreducible components of $M$ have dimension $d$, hence
$M$ is a complete intersection, \vspace{-2mm}
\item $X$ is a closed subscheme of $M$ and $X=M$ at the generic
point of every irreducible component of $X$, \vspace{-2mm}
\item some $r$-minor of the Jacobian matrix of $F_1,\ldots,F_r$
does not vanish at the generic point of any irreducible component
of $X$, \vspace{-2mm}
\end{enumerate}
and we have the same properties for any other closed subscheme
defined by $r$ elements of $\{F_1,\ldots,F_m\}$. This assumption
will be true for a general element of $GL_m(k)$.

The following Lemma is an amelioration of Lemma B and is due to
Ein and Musta\c t\v a \cite[Lemma 4.2]{EinMustata}.

\vspace{0,5cm}

\noindent \textbf{Lemma C.} Let $\rho \in \mathcal{L}(M) \setminus
\mathcal{L}(\mbox{Jac}_M)$ and put $e=e(M;\rho)$. If $\pi_{e}(\rho)
\in \mathcal{L}_{e}(X)$, then $\rho \in \mathcal{L}(X)$.

\vspace{0,5cm}

\noindent \textbf{Proposition 7.} Suppose that $n \geq
\max\{2e_r,e:=e_1+\cdots+e_r\}$. \\ (a) Let $\theta \in
\mathcal{L}^{e_1,\ldots,e_r}_n(X)$. For $l < n$, we define
\[
A_l := \{ \alpha \in \mathcal{L}_n(\mathbb{A}^N) \mid
\pi^n_l(\alpha) = \pi^n_l(\theta) \}.
\]
Take a closed subscheme $M$ associated to $r$ elements of
$\{F_1,\ldots,F_m\}$ such that $\theta \in
\mathcal{L}^{e_1,\ldots,e_r}_n(M)$. Then,
\begin{enumerate} \vspace{-3mm}
\item[(i)] we have that $\pi^n_{n-e_r}(\theta)$ is liftable
to an arc $\gamma \in \mathcal{L}(X)$,

\vspace{-3mm}
\item[(ii)] for $n > l \geq e_r$, we have
$\mathcal{L}_n(M) \cap A_l = \mathcal{L}_n(X) \cap A_l$,

\vspace{-3mm}
\item[(iii)] we have $n-e_r < e_{r+1}(X;\theta)$ if $r<N$.
\end{enumerate}
\vspace{-2mm} (b) For $n > l \geq e_r$, the set
$\pi^n_l(\mathcal{L}_n^{e_1,\ldots,e_r}(X))$ is constructible and
the truncation map \[ \mathcal{L}_n^{e_1,\ldots,e_r}(X)
\rightarrow \pi^n_l(\mathcal{L}_n^{e_1,\ldots,e_r}(X))
\]
is a piecewise trivial fibration with fibre an affine space. If $l
\leq n-e_r$, the dimension of this fibre is equal to $(n-l)d+e$.
If $l=e_r$ and $d \geq 1$, this dimension is at least $\ulcorner
(d+1)(n/2) \urcorner$.

\vspace{0,5cm}

\noindent \textsl{Remark.} An $n$-jet $\theta \in
\mathcal{L}_n(X)$ belongs to $\mathcal{L}_n(\mbox{Jac}_X)$ if and
only if $n < e(X;\theta)$.

\vspace{0,5cm}

\noindent \textbf{Proof.} (a.i) It follows from Proposition 1(a)
that $\pi^n_{n-e_r}(\theta)$ is liftable to an arc $\gamma \in
\mathcal{L}(M)$. Then, $\widetilde{\gamma} := \pi_n(\gamma) \in
\mathcal{L}_n(M)$. By Proposition 3(d), $\widetilde{\gamma} \in
\mathcal{L}_n(X)$. By using $n \geq e$ and Lemma C, we get $\gamma
\in \mathcal{L}(X)$.

(a.ii) We only have to prove that $\mathcal{L}_n(M) \cap A_{e_r}
\subset \mathcal{L}_n(X)$. Take $\alpha \in \mathcal{L}_n(M) \cap
A_{e_r}$. Then $\pi^n_{n-e_r}(\alpha)$ is liftable to an arc $\rho
\in \mathcal{L}(M)$. By the previous part $\pi^n_{n-e_r}(\theta)$
is liftable to an arc $\gamma \in \mathcal{L}(X)$. Because $e_r
\leq n-e_r$, we get $\pi_{e_r}(\rho) = \pi_{e_r}(\gamma)$, and by
using Proposition 5(a.i), we obtain $\rho \in \mathcal{L}(X)$, and
hence $\widetilde{\rho} = \pi_n(\rho) \in \mathcal{L}_n(X)$.
Because $\pi^n_{n-e_r}(\widetilde{\rho}) = \pi^n_{n-e_r}(\alpha)$,
we get $\alpha \in \mathcal{L}_n(X)$ by Proposition 3(d).

(a.iii) This follows from the previous part and from the fact that
the calculations show that we would get an extra nontrivial
condition for an element of $\mathcal{L}_n(M) \cap A_{e_r}$ to be
in $\mathcal{L}_n(X)$ if $n-e_r \geq e_{r+1}(X;\theta)$.

(b) It is clear that $\pi^n_l(\mathcal{L}_n^{e_1,\ldots,e_r}(X))$
is constructible. The other part is similar to the corresponding
part in the proof of Proposition 5(b). Now, we have to consider
the closed subschemes $M_i$ associated to $r$ elements of
$\{F_1,\ldots,F_m\}$ and use the cover
\[
\mathcal{L}_n^{e_1,\ldots,e_r}(X) = \bigcup_i \left(
\mathcal{L}_n^{e_1,\ldots,e_r}(M_i) \cap
\mathcal{L}_n^{e_1,\ldots,e_r}(X) \right)
\]
of $\mathcal{L}_n^{e_1,\ldots,e_r}(X)$, together with Proposition
2(c) and part (a.ii) of this proposition. $\qquad \Box$

\vspace{0,5cm}

\noindent \textsl{Remark.} If we consider an $n$-jet $\theta$ as a
morphism $\mbox{Spec}(k[[t]]/(t^{n+1})) \rightarrow X$, we have
that $\theta^* \Omega_X \cong M_n(X;\theta)$ as $k[[t]]$-modules,
so that we can extend the definition of $M_n(X;\theta)$ to an
arbitrary $n$-jet $\theta$ on a separated scheme of finite type
$X$. If $\theta \in \mathcal{L}_n(X)$ lies in open subsets $U$ and
$U'$ of $X$ for which we have closed immersions $U \hookrightarrow
\mathbb{A}^N$ and $U' \hookrightarrow \mathbb{A}^{N'}$ with $N'
\geq N$, then we have obviously that $e_1'=\cdots=e_{N'-N}'=0$ and
$e_i=e_{i+N'-N}'$ for $i=1,\ldots,N$. We can use this to define
the sets $\mathcal{L}_n^{e_1,\ldots,e_b}(X)$ for an arbitrary
separated scheme of finite type $X$. Proposition 6.b and 7.b imply
analogous statements in this more general context, and also the
following theorem is valid for an arbitrary separated scheme $X$
of finite type over $k$.

\vspace{0,5cm}

\noindent \textbf{Theorem 8.} Suppose that $X$ is reduced and of
pure dimension $d \geq 1$. Then $[\mathcal{L}_n(X)]$ is a multiple
of $\mathbb{L}^{\ulcorner (d+1)(n/2) \urcorner}$ in
$K_0(\mbox{Var}_{k})$ for all $n \in \mathbb{Z}_{\geq 0}$.

\vspace{0,5cm}

\noindent \textbf{Proof.} This follows from Proposition 6.b and
Proposition 7.b and the fact that in every case, at least one of
the conditions $e \geq n/2$ or $n \geq \max\{2e_r,e\}$ is
satisfied. $\qquad \Box$

\section{The smallest pole of motivic zeta functions}

From now on, we suppose that $k=\mathbb{C}$. Let $Y$ be a
nonsingular irreducible algebraic variety of dimension $\delta$.
We know that the closed subschemes of $Y$ are in one to one
correspondence with the coherent ideal sheaves on $Y$. Let $X$ be
a closed subscheme of $Y$, and let $\mathcal{I}_X$ be the
corresponding ideal sheaf. For $n \in \mathbb{Z}_{\geq 0}$, the
set
\[
\mathcal{X}_n := \{ \gamma \in \mathcal{L}_n(Y) \mid \gamma \cdot
X=n\}
\]
is a locally closed subvariety of $\mathcal{L}_n(Y)$. Note that
$\gamma \cdot X = \min \{ \mbox{ord}_t(f_n(\gamma)) \mid f \in
\mathcal{I}_X(Y) \}$ if $Y$ is an affine variety. For a regular
function $f:Y \rightarrow \mathbb{A}^1$ and $n \in
\mathbb{Z}_{\geq 0}$, we considered here the induced morphism $f_n
: \mathcal{L}_n(Y) \rightarrow \mathcal{L}_n(\mathbb{A}^1)$. The
motivic zeta function $Z(t)$ of $X$ or $\mathcal{I}_X$ is by
definition
\[
Z(t) := \sum_{n \geq 0} [\mathcal{X}_n]t^n \in
\mathcal{M}_{\mathbb{C}}[[t]],
\]
where $\mathcal{M}_{\mathbb{C}}$ is the localization of
$K_0(\mbox{Var}_{\mathbb{C}})$ with respect to $\mathbb{L}$. In
many references the motivic zeta function is defined as $\sum
[\mathcal{X}_n] \mathbb{L}^{-\delta n}t^n$, which is a
normalization of our version.

We now describe a formula for $Z(t)$ in terms of a
principalization of $\mathcal{I}_X$. Let $\sigma:W \rightarrow Y$
be a principalization of $\mathcal{I}_X$, i.e. $\sigma$ is a
proper birational morphism from a nonsingular variety $W$ such
that $\sigma$ is an isomorphism on $W \setminus \sigma^{-1}(X)$
and the total transform $\sigma^*(\mathcal{I}_X)$ is the ideal
sheaf of a simple normal crossings divisor $E$. Let $E_j$, $j\in
T$, be the irreducible components of $E$. Let $K_{W/Y}$ be the
relative canonical divisor supported on the exceptional locus of
$\sigma$. We define the numerical data $N_j$ and $\nu_j$ by the
equalities $E=\sum_{j\in T}N_jE_j$ and $K_{W/Y}=\sum_{j \in
T}(\nu_j-1)E_j$. For $J \subset T$, denote $E_J^{\circ} :=
(\cap_{j \in J} E_j) \setminus (\cup_{j \notin J} E_j)$. In terms
of these data, the announced formula for the motivic zeta function
is
\begin{eqnarray}
Z(t) = \sum_{J \subset T} [E_J^{\circ}] \prod_{j \in J}
\frac{(\mathbb{L}-1) \mathbb{L}^{\delta
N_j-\nu_j}t^{N_j}}{1-\mathbb{L}^{\delta N_j-\nu_j}t^{N_j}}.
\label{formulamotivic}
\end{eqnarray}
This formula is mentioned in \cite[Section 2.4]{VeysZuniga} and is
proved in the same way as the expression of the motivic zeta
function of a hypersurface in terms of an embedded resolution
\cite{DL1}.

Hironaka \cite{Hironaka} proved that there exists a
principalization $\sigma: W \rightarrow Y$ of $\mathcal{I}_X$
which is a composition
\[
Y=Y_0 \stackrel{\sigma_1}{\longleftarrow} Y_1
\stackrel{\sigma_2}{\longleftarrow} Y_2 \cdots
\stackrel{\sigma_i}{\longleftarrow} Y_i \cdots
\stackrel{\sigma_r}{\longleftarrow}Y_r=W
\]
of blow-ups $\sigma_i:Y_i \rightarrow Y_{i-1}$ in smooth
irreducible centers $C_{i-1} \subset Y_{i-1}$ such that
\begin{enumerate} \vspace{-2mm}
\item[(i)] if we denote the exceptional variety of $\sigma_i$ and its strict
transforms by $E_i$, we have that $\cup_{1 \leq j \leq i} E_j$ has
only simple normal crossings in $Y_i$ and $C_i$ has simple normal
crossings with $\cup_{1 \leq j \leq i} E_j$ and \vspace{-2mm}
\item[(ii)] each center $C_i$ is in the zero locus of the
weak transform of $\mathcal{I}_X$ through the map $\sigma_1 \circ
\cdots \circ \sigma_i$.
\end{enumerate}
Note that $\{1,\ldots,r\} \subset T$ and that $T$ contains one
other index for each irreducible component of $X$ of codimension
one.

For $i \in \{0,\ldots,r-1\}$, we have the formulas
$N_{i+1}=\sum_{j \in S} N_j + \mu_{C_i}$ and $\nu_{i+1} = \sum_{j
\in S} \nu_j + \delta - \mbox{card}(S)-\mbox{dim}(C_i)$, where $S$
is the set of all $j \in \{1,\ldots,i\}$ satisfying $C_i \subset
E_j$ and where $\mu_{C_i}$ is at least one as a consequence of
condition (ii) of our principalization. Actually $\mu_{C_i}$ is
the largest number $b$ satisfying $\mathcal{J}_c \subset M^b$,
where $\mathcal{J}$ is the weak transform of $\mathcal{I}_X$ in
$Y_i$, $c$ is the generic point of $C_i$ and $M$ is the maximal
ideal of the local ring $\mathcal{O}_{Y_i,c}$ (see \cite[Remark
2.7]{Villamayor}). One obtains now easily (and similar to
\cite[Proof of Theorem 2.4.0]{Segersthesis}) that $\delta
N_j-\nu_j \geq 0$ for all $j \in T$.

Although Poonen \cite{Poonen} proved that
$K_0(\mbox{Var}_{\mathbb{C}})$ is not a domain, it is still not
known whether the localization map $K_0(\mbox{Var}_{\mathbb{C}})
\rightarrow \mathcal{M}_{\mathbb{C}}$ is injective or not. Let
$I'$ be the kernel of the localization map and let
$R'=K_0(\mbox{Var}_{\mathbb{C}})/I'$, which we identify with its
image in $\mathcal{M}_{\mathbb{C}}$. Obviously, the formula
(\ref{formulamotivic}) still holds if we consider $Z(t)$ as a
power series over $R'$. Note that this is an equality in the ring
$R'[[t]]$ because $\delta N_j-\nu_j \geq 0$ for all $j \in T$.

It it not clear to us that $\cap_{k \in \mathbb{Z}_{\geq 0}}
(\mathbb{L}^k) = \{ 0 \}$ and that a number $k \in \mathbb{Z}
\setminus \{0\}$ is not a zero-divisor in $R'$. Because these
properties are needed in the proof of Lemma D, we consider the
appropriate quotient of $R'$. Let
\[
I = \bigcap_{i \in \mathbb{Z}_{\geq 0}} \{ \alpha \in R' \mid
\exists n \in \mathbb{Z} \setminus \{0\} \, : \, n \alpha \in
(\mathbb{L}^i) \}
\]
and put $R=R'/I$. The ring $R$ still specializes to Hodge-Deligne
polynomials.

From now on, we will consider the motivic zeta function $Z(t)$ as
a power series over $R$. Notice that the formula
(\ref{formulamotivic}) also holds over $R$. We write the motivic
zeta function in the form
\[
Z(t) = \frac{B(t)}{\prod_{j \in S} (1-\mathbb{L}^{\delta
N_j-\nu_j}t^{N_j})},
\]
where $S \subset T$ and where $B(t)$ is not divisible by any of
the $1-\mathbb{L}^{\delta N_j-\nu_j}t^{N_j}$, with $j \in S$. Put
$l:=\min \{-\nu_j/N_j \mid j \in S \}$. The following result is
proved in \cite[Corollary to Proposition 4.6]{SVV}.

\vspace{0,5cm}

\noindent \textbf{Lemma D.} If there exists an integer $a$ such
that $[\mathcal{X}_n]$ is a multiple of $\mathbb{L}^{\ulcorner
(\delta +l')n-a \urcorner}$ for all $n$ satisfying $(\delta
+l')n-a \geq 0$, then $l' \leq l$.

\vspace{0,5cm}

\noindent \textbf{Theorem 9.} Suppose that $Y$ is a nonsingular
irreducible algebraic variety of dimension $\delta$ and that $X$
is a closed subscheme of $Y$.

(a) If $\mathcal{I}_X$ has locally $m$ generators, with $m <
\delta$, then the motivic zeta function $Z(t) \in R[[t]]$ of $X$
belongs to the subring of $R[[t]]$ generated by $R[t]$ and the
elements $(1-\mathbb{L}^{\delta N-\nu}t^{N})^{-1}$, with $\nu,N
\in \mathbb{Z}_{>0}$ and $\nu/N \leq \delta - (\delta-m+1)/2$.

(b)  If $X$ is reduced and of pure dimension $d \geq 1$, then the
motivic zeta function $Z(t) \in R[[t]]$ of $X$ belongs to the
subring of $R[[t]]$ generated by $R[t]$ and the elements
$(1-\mathbb{L}^{\delta N-\nu}t^{N})^{-1}$, with $\nu,N \in
\mathbb{Z}_{>0}$ and $\nu/N \leq \delta - (d+1)/2$.

\vspace{0,5cm}

\noindent \textbf{Proof.} (b) We have that $[\mathcal{X}_n] =
\mathbb{L}^{\delta} [\mathcal{L}_{n-1}(X)] - [\mathcal{L}_n(X)]$
for $n \geq 1$. By using Theorem 8, we get that $[\mathcal{X}_n]$
is a multiple of $\mathbb{L}^{\ulcorner (d+1)(n/2) \urcorner}$. It
follows now from Lemma D that $l \geq -\delta + (d+1)/2$ and this
finishes the proof.

(a) This is proved in the same way as part (b). We have to use now
Theorem 4 instead of Theorem 8. Indeed, Theorem 4 implies that
$[\mathcal{L}_n(X)]$ is a multiple of $\mathbb{L}^{\ulcorner
(\delta-m+1)(n/2) \urcorner}$ in $K_0(\mbox{Var}_{\mathbb{C}})$
for all $n \in \mathbb{Z}_{\geq 0}$. $\qquad \Box$

\vspace{0,5cm}

Because $R$ specializes to Hodge-Deligne polynomials, we get a
similar result for the Hodge zeta function. We can specialize
further to the topological Euler-Poincar\'e characteristic $\chi$.
The zeta function at that level is the topological zeta function,
which is the rational function in the complex variable $s$ defined
by
\[
Z_{\mathrm{top}}(s) := \sum_{J \subset T} \chi(E_J^{\circ})
\prod_{j \in J} \frac{1}{\nu_j+sN_j}.
\]
The previous theorem implies the following.

\vspace{0,5cm}

\noindent \textbf{Theorem 10.} Suppose that $Y$ is a nonsingular
irreducible algebraic variety of dimension $\delta$ and that $X$
is a closed subscheme of $Y$.

(a) If $\mathcal{I}_X$ has locally $m$ generators, with $m <
\delta$, then the topological zeta function has no pole less than
$-\delta + (\delta-m+1)/2$.

(b) If $X$ is reduced and of pure dimension $d \geq 1$, then the
topological zeta function has no pole less than $-\delta +
(d+1)/2$.

\vspace{0,5cm}

\noindent \textsl{Example.} Let $Y=\mathbb{A}^{\delta}$ and let
$X$ be the closed subscheme of $Y$ defined by the ideal
$(x_1,\ldots,x_{r-1},x_r^2+\cdots+x_{\delta}^2)$ of
$\mathbb{C}[x_1,\ldots,x_{\delta}]$. The dimension $d$ of $X$ is
equal to $\delta-r$. One can calculate that the topological zeta
function of $X$ is
\[
Z_{\mathrm{top}}(s) = \frac{r^2+r\delta-r +
(2r-1)s}{(\delta+r-1+2s)(r+s)}
\]
if $\delta-r$ is even and
\[
Z_{\mathrm{top}}(s) = \frac{r^2+r\delta-r +
(2r-2)s}{(\delta+r-1+2s)(r+s)}
\]
if $\delta-r$ is odd. This example proves that our lower bound is
optimal.

\footnotesize{

\vspace{0,5cm}

\noindent \textsc{Universidad Complutense de Madrid, Depto. Algebra. Fac. de Ciencias Matemáticas, Plaza de las Ciencias 3, 28040 Madrid, Spain} \\
\textsl{E-mail address:} helenacobo@gmail.com

\vspace{0,5cm}

\noindent \textsc{University of Leuven, Department of
Mathematics, Celestijnenlaan 200B, B-3001 Leuven, Belgium} \\
\textsl{E-mail address:} dirk.segers@wis.kuleuven.be \\
\textsl{URL:} http://wis.kuleuven.be/algebra/segers/segers.htm
\end{document}